\documentclass[12pt,a4paper]{amsart}
\usepackage{amssymb,amsmath}
\usepackage{mathrsfs}

\headsep=1cm \numberwithin{equation}{section}
\newtheorem{theorem}{Theorem}[section]
\newtheorem{lemma}[theorem]{Lemma}
\newtheorem{proposition}[theorem]{Proposition}
\newtheorem{corollary}[theorem]{Corollary}

\theoremstyle{definition}
\newtheorem{definition}[theorem]{Definition}

\newcommand\Supp{\operatorname{Supp}}
\newcommand\Ass{\operatorname{Ass}}
\newcommand\mAss{\operatorname{mAss}}

\newcommand\Ann{\operatorname{Ann}}

\newcommand\Rad{\operatorname{Rad}}

\newcommand\height{\operatorname{height}}

\begin{document}

\title[Symbolic powers of ideals and their topology]{Symbolic powers of ideals and their topology over a module}%
\author{Adeleh Azari, Simin Mollamahmoudi and Reza Naghipour$^*$}%
\address{Department of Mathematics, University of Tabriz,
Tabriz, Iran, and School of Mathematics, Institute for Research in Fundamental
Sciences (IPM), P.O. Box: 19395-5746, Tehran, Iran.}%
\email{adeleh\_azari@yahoo.com (Adeleh Azari)}
\email{mahmoudi.simin@yahoo.com (Simin Mollamahmoudi)}
\email{naghipour@ipm.ir (Reza Naghipour)} \email {naghipour@tabrizu.ac.ir (Reza Naghipour)}

\thanks{ 
This research was in part supported by a grant from IPM. \\
$^*$Corresponding author: e-mail: {\it naghipour@ipm.ir} (Reza
Naghipour)}%

\begin{abstract}
 Let $I$ denote an ideal of a Noetherian ring $R$  and $N$ a  non-zero finitely generated $R$-module. In the present paper, some necessary and sufficient conditions are given to determine when the $I$-adic topology on $N$  is equivalent to the $I$-symbolic topology on $N$. Among other things, we shall give a complete solution to the question raised by R. Hartshorne in [{\it Affine duality and cofiniteness}, Invent. Math. {\bf9}(1970), 145-164], for a prime ideal $\frak p$ of dimension one
in a local Noetherian ring $R$, by showing that the $\frak{p}$-adic topology on $N$  is equivalent to the $\frak{p}$-symbolic topology on $N$ if and only if for all  $z\in \Ass_{R^*}N^*$ there exists $\frak{q}\in \Supp(N^*)$ such that  $z\subseteq \frak{q}$  and $\frak{q}\cap R=\frak{p}.$ Also, it is shown that if for every ${\mathfrak{p}}\in \Supp(N)$ with $\dim R/\mathfrak{p}=1$, the $\mathfrak{p}$-adic and the $\mathfrak{p}$-symbolic topologies  are equivalent on $N$, then  $N$ is unmixed and $\Ass_{R} N$ has only one element. Finally, we show that if  $\Ass_{R_{\mathfrak{p}}^*}{N^*_{\mathfrak{p}}}$  consists of a single prime ideal, for all ${\mathfrak{p}}\in {A^*}(I,N)$, then the $I$-adic and the $I$-symbolic topologies on $N$ are equivalent.
\end{abstract}
\maketitle
\section {Introduction}
In this paper we continue the study of the equivalence between  the $I$-adic topology  and the $I$-symbolic topology on $N$.    Let $R$ denote a commutative Noetherian ring, $I$ an ideal of $R$ and let $N$ be a finitely generated $R$-module. For a natural
number $n$, the $n$th symbolic power of $I$ with respect to $N,$ denoted by
$(IN)^{(n)},$ is defined to be the intersection of those primary components of
$I^nN$ which correspond to the minimal primes of
$\Ass _R N/IN.$ It is easy to see that the submodule $(IN)^{(n)}$
 is the set of all elements $x$ in $N$ for which there exists an element $s$ in
 $R\setminus\bigcup\frak{p}$ such that
 $sx\in I^nN,$ where $\frak{p}$ runs over the set of minimal primes of
$\Ass _R N/IN.$ This definition is inspired by the one given in
\cite[Remark, p. 233]{zs}. Symbolic powers of ideals are central objects in commutative
 algebra and algebraic
  geometry for their tight connection to primary decomposition of ideals and the order of vanishing
  of polynomials.

  The $I$-adic filtration   $\lbrace I^nN\rbrace_{n\geqslant0}$ and the $I$-symbolic filtration
    $\lbrace (IN)^{(n)}\rbrace_{n\geqslant0}$ induce topologies on $N$ which are called the
  $I$-adic and the $I$-symbolic topology respectively. One readily sees from the definition that
  $I^nN\subseteq(IN)^{(n)}$ for all natural numbers $n$,   so that the $I$-adic topology on $N$ is stronger than the $I$-symbolic
  topology on $N,$ but there are not equivalent in general. Therefore, one would like to compare
 the $I$-adic topology  and the $I$-symbolic topology on $N$  and provide some criterions for
 equivalence. These two topologies are said to be equivalent if for
 every natural number $n,$ there is a natural number $m$ such that  $I^nN$ contains $(IN)^{(m)}.$
R. Hartshorne in \cite[Proposition 7.1]{Ha} proved that  if $\frak p$ is a prime ideal of dimension one of a complete  local ring $R$, then the $\frak{p}$-adic topology on $N$  is equivalent to the $\frak{p}$-symbolic topology on $N$ if and only if  every associated prime ideal of $N$ is contained in $\frak p$.  In this paper Hartshorne writes: \lq \lq A general question, whose solution is quite complicated, is to determine when the $\frak{p}$-adic topology is equivalent to the $\frak{p}$-symbolic topology''. With respect to this question, P. Schenzel in \cite[Theorem 1]{sc1} gave a solution to this problem, in the case when $R=N$, and later S. McAdam and L.J. Ratliff in \cite{MR1}  gave an elegant proof of Schenzel's theorem.  Finally, L.J. Ratliff in \cite{ra1}  and J.K. Verma in \cite{Ve1} generalized Schenzel's theorem to  primary and arbitrary ideals, respectively. The purpose of this paper is to give a  generalization of Hartshorne's result by removing the completeness condition on ring. Moreover, we prove some new results concerning on the equivalence of the $I$-adic topology and the $I$-symbolic topology on $N$. Namely,  we show that these topologies are equivalent in the following  cases:

$\rm (i)$  The $\mathfrak{p}$-adic and the $\mathfrak{p}$-symbolic topologies on $N$ are equivalent, for all $\frak p\in \mAss _R N/IN$.

 $\rm (ii)$  $N$ is locally unmixed and $I=\Rad(J+\Ann_{R} N)$, where $J$ is an $N$-proper ideal of $R$ generated by  $\height_{N} J$ elements.

$\rm (iii)$  $\Ass_{R_{\mathfrak{p}}^*}{N^*_{\mathfrak{p}}}$ consists of a single prime ideal, for all ${\mathfrak{p}}\in {A^*}(I,N)$.

$\rm (iv)$  For every ideal $J$ of $R$ with $\Ass _R N/JN=\mAss _R N/JN$ and $I\subseteq J$, the $I$-symbolic topology is finer than the $J$-symbolic topology on $N$.

In addition, we show that if $R$ is local and for every ${\mathfrak{p}}\in \Supp(N)$ with $\dim R/\mathfrak{p}=1$, the $\mathfrak{p}$-adic and the $\mathfrak{p}$-symbolic topologies  are equivalent, then  $N$ is unmixed and $\Ass_{R} N$ has only one element.\\

 Throughout this paper all rings are commutative and Noetherian, with
 identity, unless otherwise specified. We shall use $R$ to denote
 such a ring, $I$ an ideal of $R,$ and $N$ a non-zero finitely
 generated module over $R.$ We denote by
 $\mathscr{R}$  the {\it Rees ring}
 $R[u,It]:=\oplus_{n\in\mathbb{Z}}I^nt^n$ of $R$ with respect to $I,$ where $t$ is an indeterminate and
 $u=t^{-1}.$ Also, the {\it graded Rees module}
 $N[u,It]:=\oplus_{n\in\mathbb{Z}}I^nN$ over $\mathscr{R}$ is denoted by
 $\mathscr{N},$ which is a finitely  generated $\mathscr{R}$-module.

 If $(R,\frak{m})$ is local, then $R^{*}$ (resp. $N^{*}$)
 denotes the completion of $R$ (resp. $N$) w.r.t. the $\frak{m}$-adic topology.
Then $N$ is said to be an {\it unmixed}  (resp. a {\it quasi-unmixed}) {\it module} if
all the prime ideals  (resp. all the minimal prime ideals) of
$\Ass_{R^*}N^*$ have the same dimension. More generally, if $R$ is
not necessarily local, $N$ is a {\it locally unmixed} (resp. {\it locally
quasi-unmixed}) {\it module}  if for any $\frak{p}\in\Supp(N),$
$N_{\frak{p}}$ is an unmixed (a quasi-unmixed) module over
$R_\frak{p}.$ We shall say that an ideal $J$ of $R$ is $N$-{\it proper} if
$N/JN\neq0,$ and, when this is the case, we define the $N$-$\height$ of
$J$ (written $\height_N J $) to be $\inf \lbrace \height_N\frak{p}:
\frak{p}\in\Supp{N}\cap V(J)\rbrace,$ where $ \height_N\frak{p}$ is
defined to be the supremum of lengths of chains of prime ideals of
$\Supp(N)$ terminating with $\frak{p}.$ For any ideal $\frak{a}$ of
$R,$ the {\it radical} of $\frak{a},$ denoted by $\Rad(\frak{a}),$
 is defined to be the set
 $\lbrace x\in\frak{a}: x^n\in \frak{a}\,\,  \text {for some}\,  n\in \mathbb{N}\rbrace.$ For any
 $R$-module $L,$ we denote by $\mAss_RL$ the set of minimal prime ideals of $\Ass_RL.$ For any
unexplained notation and terminology we refer the reader to \cite{Ma} or \cite{Nag}.\\

 In the second section, we examine the equivalence  between the $I$-adic and the $I$-symbolic topologies.  In this section, we show that these topologies are equivalent, whenever for every
 $\frak{p}\in \mAss _R N/IN,$ the $\frak{p}$-adic topology on $N$  is equivalent to
 the $\frak{p}$-symbolic topology on $N.$

 The main result of the third section is that if $N$ is a non-zero finitely generated locally unmixed
 $R$-module and $I$ an $N$-proper ideal of $R$ generated by $\height_N I$ elements, then
 the $J$-adic topology is equivalent to
 the $J$-symbolic topology, where $J=\Rad(I+\Ann_RN).$ Also, in this section we shall extend a nice result of  Hartshorne \cite[Proposition 7.1]{Ha}.

 Namely,  we shall show the following result:

\begin{theorem}
 If  $(R,\frak{m})$ is local and $\frak{p}\in \Supp(N)$ with $\dim R/\frak{p}=1,$ then the $\frak{p}$-adic topology on $N$  is equivalent to the $\frak{p}$-symbolic topology on $N$ if and only if for all
 $z\in \Ass_{R^*}N^*$ there exists $\frak{q}\in \Supp(N^*)$ such that $z\subseteq \frak{q}$  and $\frak{q}\cap R=\frak{p}.$
\end{theorem}

 Finally, the main purpose of the Section $4$ is to establish a connection between the  unmixedness (resp. associated primes) of $N$ and the comparison of the
 $\frak p$-adic topology on $N$  and the $\frak p$-symbolic topology on $N$  (resp.  $I$-adic topology and the $I$-symbolic topology on $N$),  for certain prime ideals $\frak p$ of $R$.  More precisely, we show that:

 \begin{theorem}
Let  $N$ be a non-zero finitely generated module over a commutative Noetherian ring $R$.
\begin{itemize}
  \item[(i)] If $\Ass_{R^*_\frak{p}}N^*_\frak{p}$ consists of a single prime, for all
 $\frak{p}\in A^*(I),$ then the $I$-adic topology on $N$  is equivalent to the $I$-symbolic topology on $N$.

\item[(ii)] If $(R,\frak{m})$ is local and  for every $\frak{p}\in \Supp(N)$ with $\dim R/\frak{p}=1,$  the $\frak{p}$-adic topology on $N$ is equivalent to
 the $\frak{p}$-symbolic topology on $N,$ then $N$ is unmixed and $\Ass_R N$ consists of a single prime.
 \end{itemize}
 \end{theorem}
\section{Comparison of Topologies}
The purpose of this section is to examine the equivalence  between the $I$-adic and the $I$-symbolic topologies.  The main goal of this
section is Theorem 2.5, which shows that these topologies are equivalent, whenever for every  $\frak{p}\in \mAss _R N/IN,$ the
$\frak{p}$-adic topology is equivalent to  the $\frak{p}$-symbolic topology. Before we state the main result of this section, let us give a definition:
\begin{definition} \label{1.1}
A prime ideal $\frak{p}$ of $R$ is called a {\it quitessential prime ideal of} $I$
w.r.t. $N$ precisely when there exists $\frak{q}\in\Ass_{R^*_\frak{p}}N^*_\frak{p}$ such that
$\Rad(IR^*_\frak{p}+ \frak{q})= \frak{p}R^*_\frak{p}$. The set of quitessential  primes of $I$
w.r.t. $N$ is denoted by $Q(I,N)$ which is a finite set.  Then the {\it essential primes of }$I$ w.r.t. $N$, denoted by $E(I,N)$, is defined to be the
set $\{\frak{q}\cap R\mid\, \frak{q}\in Q((u\mathscr{R},\mathscr{N}))\}$.
\end{definition}
The concepts of the quitessential and essential prime ideals of $I$ were introduced by McAdam \cite{Mc2}, and  Ahn  in \cite{Ah} extended they to modules.

\begin{proposition}\label{prop1}
  Let  $N$ be a finitely generated $R$-module and $I$  an $N$-proper ideal of $R$. Then the following conditions are equivalent:

\begin{itemize}
  \item[(i)]  The $I$-symbolic topology on $N$ is  equivalent  to the  $I$-adic topology on $N$.

\item[(ii)]   For every $N$-proper ideal $J$ of $R$ with $\Ass _R N/JN=\mAss _R N/JN$ and $I\subseteq J$,
the $I$-symbolic  topology on $N$  is finer than the $J$-symbolic  topology on $N$.
\end{itemize}
\end{proposition}
\proof $\rm(i) \Longrightarrow (ii)$:  Let $J$ be an $N$-proper ideal of $R$ such that  $I\subseteq J$ and that $\Ass _R N/JN=\mAss _R N/JN$. For every  integer $l\geq0$, we need to show that there exists an integer $m\geq0$ such that $(IN)^{(m)}\subseteq(JN)^{(l)}$. To this end, in view of the assumption (i), there is an integer $m\geq0$ such that
 $(IN)^{(m)}\subseteq I^lN$, and so $(IN)^{(m)}\subseteq J^lN\subseteq(JN)^{(l)}$, as required.

 $\rm(ii) \Longrightarrow (i)$:  In view of \cite[Corollary 3.7]{Na1}, it is enough for us to show that $Q(I,N)=\mAss _R N/IN$. To achieve this, suppose the contrary is true. Then there is an element ${\mathfrak{p}}\in Q(I,N)$ such that ${\mathfrak{p}}\notin \mAss _R N/IN$. (Note that $\mAss _R N/IN\subseteq Q(I,N)$.) Hence, in view of  \cite[Theorem 3.6]{Na1}, there exists an integer $k\geq0$ such that $(IN)^{(m)}\nsubseteq(\mathfrak{p}N)^{(k)}$ for all integer $m\geq 0$. Now, because of $I\subseteq {\mathfrak{p}}$ and $\Ass _R N/{\mathfrak{p}}N=\mAss _R N/{\mathfrak{p}}N$,  the assumption (ii) provides a contradiction. \qed\\

The following proposition and its corollary  are quite useful in the proof of the main theorem.

\begin{proposition}\label{prop2}
  Let   $N$ be a non-zero finitely generated $R$-module and $I_1, I_2$ be two $N$-proper ideals of
$R$ such that
\begin{equation*}\label{eq1}
  \mAss _R N/(I_1I_2)N=\mAss _R N/I_1N \cup \mAss _R N/I_2N .
\end{equation*}
Suppose that the $I_i$-symbolic topology on $N$ is equivalent to the $I_i$-adic topology on $N$,  for $i=1, 2$.
Then the $I_1I_2$-symbolic topology on $N$ is equivalent to the  $I_1I_2$-adic topology on $N$ and  $(I_1\cap I_2)$-symbolic topology on $N$ is equivalent to the $(I_1\cap I_2)$-adic topology  on $N$.
\end{proposition}
\proof As $\Rad(I_1I_2)=\Rad(I_1\cap I_2)$ and $\mAss _R N/(I_1I_2)N=\mAss _R N/(I_1\cap I_2)N$, it is enough, in view of \cite[Lemma 3.1 and Corollary 3.7]{Na1}, to show that $Q(I_1I_2,N)=\mAss _R N/(I_1I_2)N$. To achieve this, suppose that  ${\mathfrak{p}}\in Q(I_1I_2,N)$. Then there exists $z\in\Ass_{R_{\mathfrak{p}}^*}{N^*_{\mathfrak{p}}}$ such that $\Rad(I_1I_2{R_{\mathfrak{p}}^*}+z)={\mathfrak{p}}R_{\mathfrak{p}}^*$. As $I_1I_2\subseteq {\mathfrak{p}}$, without loss of generality we may assume that $I_1\subseteq {\mathfrak{p}}$. Then $\Rad(I_1{R_{\mathfrak{p}}^*}+z)={\mathfrak{p}}R_{\mathfrak{p}}^*$, and so ${\mathfrak{p}}\in Q(I_1,N)$. Hence, in view of assumption and \cite[Corollary 3.7]{Na1},   ${\mathfrak{p}}\in \mAss _R N/I_1N$. Therefore, it follows from
$$ \mAss _R N/(I_1I_2)N=\mAss _R N/I_1N \cup \mAss _R N/I_2N,$$
that ${\mathfrak{p}}\in \mAss _R N/(I_1I_2)N$, as required. \qed \\

\begin{corollary}\label{coro1}
  Let   $N$ be a non-zero  finitely generated $R$-module and let $I_1,\ldots, I_n$ be $N$-proper ideals of
$R$ such that $\mAss _R N/(\prod\limits_{i=1}^n I_i)N=\bigcup\limits_{i=1}^n \mAss _R N/I_iN$, and
that the $I_i$-symbolic topology on $N$ is equivalent to  the $I_i$-adic topology on $N$,  for all $i=1,\ldots, n$. Then the $\prod\limits_{i=1}^n I_i$-symbolic topology on $N$  (resp. $\bigcap \limits_{i=1}^n I_i$-symbolic topology on $N$) is equivalent to  the $\prod\limits_{i=1}^n I_i$-adic (resp. $\bigcap \limits_{i=1}^n I_i$-adic) topology on  $N$.
\end{corollary}
\proof The result follows from Proposition \ref{prop2} and induction on $n$. \qed \\

We are now ready to state and prove the main theorem of this section, which gives us a criterion of the equivalence  between the $I$-adic and the $I$-symbolic topologies.
\begin{theorem}\label{coro2}
Let $N$ be a non-zero finitely generated $R$-module and let $I$ be an $N$-proper ideal of
$R$ such that  the  $\mathfrak{p}$-symbolic topology  is equivalent to the $\mathfrak{p}$-adic topology on $N$,   for all
${\mathfrak{p}}\in \mAss _R N/IN$. Then the $I$-symbolic topology on $N$  is equivalent to the $I$-adic topology on $N$.
\end{theorem}
\proof Let $\mAss _R N/IN=\{\mathfrak{p}_1,\ldots, \mathfrak{p}_n\}$. Then it is easy to see that
$$\mAss _R N/(\prod\limits_{i=1}^n \mathfrak{p}_i)N=\bigcup\limits_{i=1}^n \mAss _R N/\mathfrak{p}_iN .$$
Now, the assertion follows from Corollary \ref{coro1} and  \cite[Lemma 3.1]{Na2}. \qed \\

\section{Locally Unmixed Modules and Comparison of Topologies}

The main goal of this  section is  to prove the equivalence between the $I$-adic and the $I$-symbolic topologies on a  finitely generated locally unmixed
 $R$-module $N$ for certain ideals $I$  of $R$. Also, we explore an  equivalence  between the $\frak p$-adic and the $\frak p$-symbolic topologies and the associated primes of $N$, for prime ideals $\frak p$ of dimension one.  We begin with:

 \begin{definition} \label{1.1}
 A prime ideal $\frak{p}$ of $R$ is called a  {\it quitasymptotic prime ideal of }$I$ w.r.t. $N$  precisely when there exists
$\frak{q}\in \mAss_{R^*_\frak{p}}N^*_\frak{p}$ such that $\Rad(IR^*_\frak{p}+ \frak{q})= \frak{p}R^*_\frak{p}$. The set of
quitasymptotic prime ideals of $I$ w.r.t. $N$ is denoted by $\bar{Q}(I, N)$. Then
the {\it asymptotic prime ideals of} $I$ w.r.t. $N,$ denoted by $\bar{A^{*}}(I,N)$, is defined to be the set
$\{\frak{q}\cap R\mid\, \frak{q}\in \bar{Q^{*}}(u\mathscr{R},\mathscr{N}) \}$.
\end{definition}
\begin{lemma}\label{lem1}
Let $N$ be a non-zero  finitely generated  locally unmixed $R$-module and  let $I$ be an ideal of $R$. Then $E(I,N)=\bar{A}^*(I,N)$.
\end{lemma}
\proof In view of \cite[Corollary 3.7]{Ah}, it is enough for us to show that $E(I,N)\subseteq \bar{A}^*(I,N)$. To do this let ${\mathfrak{p}}\in E(I,N)$. Since both $E(I,N)$ and $\bar{A}^*(I,N)$ behave  well under localization, without loss of generality, we may assume that $(R,\mathfrak{p})$ is local.  Also, in view \cite[Proposition 3.8]{Ah}, it is easy to see that we may assume in addition that $R$  is complete. Now, according to \cite[Proposition 3.6]{Ah}, there exists $z\in \Ass _R N$  such that $z\subseteq {\mathfrak{p}}$ and $\mathfrak{p}/z\in E(\mathfrak{p}+z/z)$. Since $R/z$ is unmixed it follows from  \cite[Proposition 2.11]{KR}  that  $\mathfrak{p}/z\in \bar{A}^*(I+z/z)$. Moreover, since by hypothesis $z\in \mAss _R N$,  it follows from  \cite[Proposition  3.6]{Ah} that  ${\mathfrak{p}}\in \bar{A}^*(I,N)$, as required. \qed \\

Before we state the next lemma, let us recall the important notion {\it analytic spread of $I$ with respect to $N$}, over a local ring $(R, \frak{m}$), introduced by Brodmann in \cite{Br2}:
$$l(I,N):= \dim\,\mathcal{N}(I,N)/ (\frak{m},u)\mathcal{N}(I,N),$$ in the case $N=R$, $l(I,N)$ is the classical analytic spread $l(I)$
of $I$, introduced by Northcott and Rees (see \cite{NR}).
\begin{lemma}\label{lem2}
Let $N$ be a non-zero  finitely generated  locally quasi-unmixed $R$-module and let $I$ be an
$N$-proper ideal of $R$ generated by $\height _NI$ elements. Then
$\bar{A}^*(I,N)=\mAss _R N/IN$.
\end{lemma}
\proof As $\mAss _R N/IN\subseteq \bar{A}^*(I,N)$, it  will suffice for us to show that $\bar{A}^*(I,N)\subseteq \mAss _R N/IN$. To this end, let ${\mathfrak{p}}\in \bar{A}^*(I,N)$. Since $N_{\mathfrak{p}}$ is a quasi-unmixed $R_{\mathfrak{p}}$-module, it follows from \cite[Proposition  2.3]{NS} that $\height _N\mathfrak{p}=\ell(IR_{\mathfrak{p}},N_{\mathfrak{p}})$. Since at least $\ell(\frak{a})$ elements are needed to generated $\frak{a}$, for any ideal $\frak{a}$ in a commutative Noetherian ring $A$, it follows from \cite[Lemma 2.2]{NS} that $\ell(IR_{\mathfrak{p}},N_{\mathfrak{p}})\leqslant \height _NI$, and so $\height _N\mathfrak{p}=\height _NI$. Therefore ${\mathfrak{p}}\in \mAss _R N/IN$, as required. \qed \\

The following theorem, which is one of our main results of this section, shows that for certain ideals $I$, the $I$-symbolic topology on $N$ is equivalent to the $I$-adic topology on  $N$, whenever $N$ is a  finitely generated  locally unmixed $R$-module.

\begin{theorem}\label{thm1}
Let $N$ be a non-zero  finitely generated  locally unmixed $R$-module and  let $J$ be an $N$-proper ideal of $R$ generated by $\height _NJ$ elements.
 Then the $I$-symbolic topology on $N$ is equivalent to the $I$-adic topology on $N$, where $I=\Rad(J+\Ann _RN)$.
\end{theorem}
\proof
In view of \cite[Corollary 3.7 and Lemma 3.1]{Na1}, it will suffice  to show that $Q(I,N)\subseteq \mAss _R N/IN$. For this let $ {\mathfrak{p}}\in Q(I,N)$. Then, it follows from  \cite[Lemma 3.1]{Na1} that  $ {\mathfrak{p}}\in Q(J,N)$. Thus, by \cite[Corollary 3.7]{Ah}, $ {\mathfrak{p}}\in E(J,N)$, and so by virtue of Lemma \ref{lem1}, ${\mathfrak{p}}\in\bar{A}^*(J,N)$. Therefore, in view of Lemma \ref{lem2}, ${\mathfrak{p}}\in \mAss _R N/JN$. Now, as
$\mAss _R N/IN=\mAss _R N/JN$,  the desired result follows. \qed \\

The next theorem, which is the final main result of this section, extends a nice result of  Hartshorne \cite[Proposition 7.1]{Ha}.
 \begin{theorem}\label{thm2}
Let $(R,\mathfrak{m})$ be a  local ring and  $N$  a non-zero  finitely generated  $R$-module. Let ${\mathfrak{p}}\in \Supp(N)$
with $\dim R/\mathfrak{p}=1$.  Then the following conditions are equivalent:
\begin{itemize}
  \item[(i)] The $\mathfrak{p}$-symbolic topology on $N$  is equivalent to the $\mathfrak{p}$-adic topology on $N$.
  \item[(ii)] For all $z\in\Ass _{R^*} N^*$ there exists $\mathfrak{q}\in \Supp(N^*)$ such that $z\subseteq\mathfrak{q}$ and $\mathfrak{q}\cap R=\mathfrak{p}$.
\end{itemize}
\end{theorem}
\proof $\rm(i) \Longrightarrow (ii)$:  Let $z\in\Ass _{R^*} N^*$. In view  of \cite[Corollary 3.7]{Na1} and the assumption (i), we  have $$Q({\mathfrak{p}},N)=\mAss _R N/{\mathfrak{p}}N=\{{\mathfrak{p}}\}.$$ Therefore  $\mathfrak{m}\notin Q({\mathfrak{p}},N)$, and so
$\mathfrak{m}R^*$ is not minimal over ${\mathfrak{p}}R^*+z$. Let $\mathfrak{q}$ be a minimal prime over ${\mathfrak{p}}R^*+z$. Then  $\mathfrak{q}\in \Supp(N^*)$ and ${\mathfrak{p}}\subseteq \mathfrak{q}\cap R$. Now, as $\dim R/\mathfrak{p}=1$, one easily sees  that  $\mathfrak{q}\cap R=\mathfrak{p}$ and $z\subseteq \mathfrak{q}$, as required.

In order to prove $\rm(ii) \Longrightarrow (i)$, in view  of \cite[Corollary 3.7]{Na1}, it is enough for us to show that  $Q({\mathfrak{p}},N)=\{\mathfrak{p}\}$. To do this, let $\mathfrak{q}\in Q({\mathfrak{p}},N)$. Then ${\mathfrak{p}}\subseteq\mathfrak{q} \subseteq \mathfrak{m}$. Since $\dim R/\mathfrak{p}=1$, we see that $\mathfrak{q}=\mathfrak{p}$ or $\mathfrak{q}=\mathfrak{m}$. If $\mathfrak{q}=\mathfrak{m}$, then $\mathfrak{m}\in Q({\mathfrak{p}},N)$, and so there exists $z\in\Ass _{R^*} N^*$ such that $\Rad({\mathfrak{p}} R^*+z)=\mathfrak{m}R^*$. Hence, by the assumption (ii), there exists  $\mathfrak{q'}\in \Supp(N^*)$ such that $z\subseteq {\mathfrak{q'}}$ and $\mathfrak{q'}\cap R=\mathfrak{p}$. Therefore, $\mathfrak{q'}\subseteq {\mathfrak{p}}R^*$, and so $\Rad({\mathfrak{p}} R^*)=\mathfrak{m}R^*$. Consequently, $\dim R^*/{\mathfrak{p}}R^*=0$ which is a contradiction, because $\dim R^*/{\mathfrak{p}}R^*=\dim R/\mathfrak{p}=1$. Whence $\mathfrak{q}=\mathfrak{p}$ and this completes the proof.  \qed\\
\section{Associated primes and Unmixedness}

The main aim of this section shows that if  $\Ass_{R^*_\frak{p}}N^*_\frak{p}$ consists of a single prime, for all
 $\frak{p}\in A^*(I, N),$ then the $I$-adic topology is equivalent to the $I$-symbolic topology on $N$.
Furthermore, we show that, if $(R,\frak{m})$ is local and  for every $\frak{p}\in \Supp(N)$ with $\dim R/\frak{p}=1,$
the $\frak{p}$-adic topology is equivalent to  the $\frak{p}$-symbolic topology on $N$, then $N$ is unmixed and $\Ass_R N$ has only one element.
Following  \cite{Br}, we shall use $A^*(I, N)$ to denote the ultimately constant values of $\Ass_R N/I^nN$ for all large $n$. The following theorem is the first main result of this section.

\begin{theorem}\label{thm3}
Let  $N$ be a non-zero  finitely generated  $R$-module and $I$  an $N$-proper ideal of $R$  such that $\Ass _{R^*_{\mathfrak{p}}} N^*_{\mathfrak{p}}$ consists of a single prime ideal $z$, for all ${\mathfrak{p}}\in A^*(I,N)$. Then the $I$-symbolic topology on $N$ is equivalent to the $I$-adic topology on $N$.
\end{theorem}
\proof In view of \cite[Corollary 3.7]{Na1}, it will suffice to show that  $Q(I,N)=\mAss _R N/IN$. To do this, suppose the contrary is true. That is there exists ${\mathfrak{p}}\in Q(I,N)$ such  that ${\mathfrak{p}}\notin \mAss _R N/IN$. Since  ${\mathfrak{p}}\in \Supp(N/IN)$, it follows that there exists $\mathfrak{q}\in \mAss _R N/IN$ such that $\mathfrak{q}\subsetneqq\mathfrak{p}$. Moreover, by virtue of \cite[Theorem 3.17]{Ah}, $Q(I,N)\subseteq A^*(I,N)$, hence $\Ass _{R^*_{\mathfrak{p}}} N^*_{\mathfrak{p}}=\{z\}$. Therefore, $\Rad(IR^*_{\mathfrak{p}}+z)={\mathfrak{p}}R^*_{\mathfrak{p}}$. Now, let $\mathfrak{q}^*$ be a minimal prime over $\mathfrak{q}R_{\mathfrak{p}}^*$.
Then  $IR^*_{\mathfrak{p}}\subseteq\mathfrak{q}R^*_{\mathfrak{p}}\subseteq\mathfrak{q}^*$. Furthermore, as $\mathfrak{q}\in \Supp(N)$, it easily follows from \cite[Theorem 18.1]{Nag} that  $\mathfrak{q}^*\in \Supp(N^*_{\mathfrak{p}})$, and so $z\subseteq \mathfrak{q}^*$. Consequently ${\mathfrak{p}}R^*_{\mathfrak{p}}\subseteq\mathfrak{q}^*$,  and hence ${\mathfrak{p}}R^*_{\mathfrak{p}}\subseteq\mathfrak{q}^*\cap R_{\mathfrak{p}}$.
On the other hand, since $\mathfrak{q}^*$ is a minimal prime over $\mathfrak{q}R^*_{\mathfrak{p}}$, we can therefore deduce from the Going-down Theorem (see \cite[Theorem 9.5]{Ma}) that $\mathfrak{q}^*\cap R_{\mathfrak{p}}=\mathfrak{q}R_{\mathfrak{p}}$. Hence $ \mathfrak{q}R_{\mathfrak{p}}={\mathfrak{p}}R_{\mathfrak{p}}$, and so $\mathfrak{q}=\mathfrak{p}$, which is a contradiction.  \qed \\

The following proposition is needed in the proof of the second main theorem.
\begin{proposition}\label{prop3}
Let $(R,\mathfrak{m})$ be a  local ring and  $N$  a non-zero  finitely generated  $R$-module such that $\dim N>0$ and that $\Ass_R N$ has at least two
 elements. Then  there exists an $N$-proper ideal $I$ of $R$ such that
\begin{eqnarray*}
\mathfrak{m}\in Q(I,N)\backslash \mAss_RN/IN .
\end{eqnarray*}
\end{proposition}
\proof In view of assumption, there exist ${z}_1, z_2\in \Ass_RN$ such that $z_1\neq z_2$. Without loss of generality, we may assume that $z_1\in \mAss_RN$. Let $n:=\dim R/z_1+z_2$. If $n=\dim N$, then there exists a minimal prime $\mathfrak{p}$ over $z_1+z_2$ such that $$\dim R/z_1+z_2=\dim R/\mathfrak{p}=n=\dim N.$$ Hence ${\mathfrak{p}}\in \mAss_RN$ and $z_1+z_2\subseteq {\mathfrak{p}}$. Consequently, $z_1=\mathfrak{p}=z_2$, which is a contradiction. Therefore, $n<\dim N$. Suppose $n=0$. Then $z_1+z_2$ is $\mathfrak{m}$-primary, and so in view of \cite[Lemma 3.5]{Ah}, $\mathfrak{m}\in Q(z_1,N)\backslash \mAss_RN/z_1N$, as required.

Now, suppose $n>0$. Then there exist  elements $a_1,\ldots,a_n$ of $\mathfrak{m}$ such that their images in $R/z_1+z_2$ form a system of parameters. Let $J=(a_1,\ldots,a_n)$. Then $\Rad(J+z_1+z_2)=\mathfrak{m}$, and so $J+z_1+z_2$ is $\mathfrak{m}$-primary. Now, if $\Rad(J+z_1)=\mathfrak{m}$, then as $z_1\in \Ass_RN$, it follows  from \cite[Lemma 3.5]{Ah} that $\mathfrak{m}\in Q(I,N)$. Moreover, $\mathfrak{m}\notin \mAss_RN/JN$. Because, if $\frak m\in \mAss_RN/JN$, then $\mathfrak{m}=\Rad(J+\Ann_RN)$. Hence, $$\height(\mathfrak{m}/\Ann N)=\height (J+\Ann_RN/\Ann_R N),$$ and so $\dim N=\height_N J\leq n$, which is a contradiction. Also, if  $\Rad(J+z_1)\neq\mathfrak{m}$, then $\mathfrak{m}\notin \mAss_R N/(J+z_1)N$. Hence, using   \cite[Lemma 3.5]{Ah} and  $\Rad(J+z_1+z_2)=\mathfrak{m}$, we obtain that $\frak m\in Q(J+z_1,N)$. This completes the proof.  \qed \\

Now, we can state and prove the second main theorem of this section.

\begin{theorem}\label{thm4}
Let $(R,\mathfrak{m})$ be a  local ring and let $N$  a non-zero  finitely generated  $R$-module of positive dimension.
Suppose that  for every ${\mathfrak{p}}\in \Supp(N)$ with  $\dim R/\mathfrak{p}=1$, the $\mathfrak{p}$-symbolic topology on $N$  is equivalent to the $\mathfrak{p}$-adic topology on $N$. Then  $N$ is unmixed and $\Ass_R N$ has exactly one element.
\end{theorem}
\proof If $\mathfrak{m}R^*\in \Ass_{R^*}N^*,$ then it  follows  easily from \cite[Theorem 23.2]{Ma} that $\mathfrak{m}\in \Ass_RN$, and  so $\mathfrak{m}\in Q(I,N)$ for  every ideal $I$ of $R$. Now, as $\dim N>0$ there exists ${\mathfrak{p}}\in \Supp(N)$ such that
$\mathfrak{p}\subsetneqq \frak m$ and  $\dim R/\mathfrak{p}=1$.
Hence $$\mathfrak{m}\in Q({\mathfrak{p}},N)=\mAss_RN/{\mathfrak{p}}N=\{{\mathfrak{p}}\},$$  which is a contradiction. Therefore $\mathfrak{m}R^*\notin \Ass_{R^*}N^*$. Now, we show that $N$ is unmixed. To do this, suppose the contrary, i.e., $N$ is not unmixed. Then there exists $z\in \Ass_{R^*}N^*$ such that  $\dim R^*/z<\dim N$. Since $\mathfrak{m}R^*\notin \Ass_{R^*}N^*$, we have $\dim R^*/z>0$. Therefore, in view of \cite[Proposition 3.5]{Na2} there exists an $N$-proper ideal $I$ of $R$ generated by $\height_NI$ elements, such that
\begin{eqnarray*}
\Rad(IR^*+z)=\mathfrak{m}R^*\quad \textrm{and}\quad \height_NI=\dim R^*/z.
\end{eqnarray*}
Consequently, $\mathfrak{m}\in Q(I,N)$, and as $\height_NI<\dim N$, there exists  ${\mathfrak{p}}\in \Supp(N)$ such that $I\subseteq {\mathfrak{p}}$ and $\dim R/\mathfrak{p}=1$. Hence $\mathfrak{p}\subsetneqq \mathfrak{m}$ and $\mathfrak{m}\in Q({\mathfrak{p}},N)$. Now, since
 $$Q({\mathfrak{p}},N)=\mAss_RN/{\mathfrak{p}}N=\{{\mathfrak{p}}\},$$ it follows that $\mathfrak{m}=\mathfrak{p}$, which is a contradiction, so $N$ is unmixed.
Now in order to complete the proof,  we must show that $\Ass_RN$ consists of a single prime. To  this end, suppose that the contrary is true.  Then, by Proposition \ref{prop3}, there exists an
ideal $I$ of $R$ such that
\begin{eqnarray*}
\mathfrak{m}\in Q(I,N)\backslash \mAss_R N/IN .
\end{eqnarray*}
Since $\mathfrak{m}\notin \mAss_R N/IN,$ there exists ${\mathfrak{p}}\in \Ass_R N/IN$ such that $\dim R/\mathfrak{p}=1$. Hence, $\mathfrak{m}\in Q({\mathfrak{p}},N)$.  Now, because of  $$Q({\mathfrak{p}},N)=\mAss_R N/{\mathfrak{p}}N=\{{\mathfrak{p}}\},$$ we see that $\mathfrak{m}=\mathfrak{p}$. Therefore, we have arrived at a contradiction, and so $\Ass_R N$ has only one element, as required.  \qed \\

\begin{center}
{\bf Acknowledgments}
\end{center}
The authors would like to thank Professor Monireh Sedghi for reading of the first draft and valuable discussions.
Finally, the authors would like to thank from the Institute for Research in Fundamental Sciences (IPM), for the financial support.


\end{document}